\documentclass[a4paper,fleqn]{cas-dc}

\usepackage{tikz}
\usepackage{tabularx}
\usepackage{float}

\let\printorcid\relax

\usetikzlibrary{decorations.markings} % <-- you need this for markings
\tikzset{->-/.style={decoration={
  markings,
  mark=at position #1 with {\arrow{>}}},postaction={decorate}}}

\usepackage[numbers]{natbib}
\def\tsc#1{\csdef{#1}{\textsc{\lowercase{#1}}\xspace}}
\tsc{WGM}
\tsc{QE}
\begin{document}
\let\WriteBookmarks\relax
\def\floatpagepagefraction{1}
\def\textpagefraction{.001}

\shorttitle{Modelling an electrolyser in a graph-based framework}    

\shortauthors{Preprint}  

\title [mode = title]{Modelling an electrolyser in a graph-based framework}  




\author[1]{Buu-Van Nguyen}[orcid=0000-0002-2810-064X]

\cormark[1]


\ead{b.nguyen@tudelft.nl}


\credit{Conceptualization, Methodology, Software, Writing - Original Draft}

\affiliation[1]{organization={Delft Insitute of Applied Mathematics, Technische Universiteit Delft},
                city={Delft},
                % citysep={}, % Uncomment if no comma needed between city and postcode
                % state={Zuid-Holland},
                country={The Netherlands}}

\author[2]{Johan Romate}%[]

% Footnote of the second author

% Email id of the second author
% \ead{}


\credit{Writing - review \& editing, Supervision}

\affiliation[2]{organization={Technische Universiteit Delft},
                city={Delft},
                % citysep={}, % Uncomment if no comma needed between city and postcode
                % state={Zuid-Holland},
                country={The Netherlands}}

\author[3]{Cornelis Vuik}[orcid=0000-0002-5988-9153]

% Footnote of the second author
%\fnmark[3]

% Email id of the second author
% \ead{c.vuik@tudelft.nl}


\credit{Writing - review \& editing, Supervision}

\affiliation[3]{organization={Delft Insitute of Applied Mathematics, Technische Universiteit Delft},
                city={Delft},
                % citysep={}, % Uncomment if no comma needed between city and postcode
                % state={Zuid-Holland},
                country={The Netherlands}}

\cortext[1]{Corresponding author}


\nonumnote{}

% Here goes the abstracft
\begin{abstract}
We propose an electrolyser model for steady-state load flow analysis of multi-carrier energy networks, where the electrolyser is capable of producing hydrogen gas and heat. We show that there are boundary conditions that lead to a well-posed problem. We derive these conditions for two cases, namely with a fixed and non-fixed ratio between gas and heat output. Furthermore, the derived conditions are validated numerically.  \nocite{*}%% Remove this line from your manuscript.
\end{abstract}


% \begin{highlights}
% \item An electrolyser model for steady-state load flow analysis
% \item Conditions for well-posedness of the electrolyser model
% \item Steady-state load flow analysis of a multi-carrier energy system including an electrolyser
% \end{highlights}

%\nocite{*}

% Keywords

\maketitle

% Main text

%%%%%%%%%%%%%%%%%%%%%%%%%%%%%%%%%%%%%%%%%%%%%%%%%%%%%%%%%%%%%%%%%%%%%%%%%%%%%%%%%%%%%%%%%%%%%%%%%%%%

% Introduction

%%%%%%%%%%%%%%%%%%%%%%%%%%%%%%%%%%%%%%%%%%%%%%%%%%%%%%%%%%%%%%%%%%%%%%%%%%%%%%%%%%%%%%%%%%%%%%%%%%%%

\section{Introduction}
\label{section:introduction}
The current level of greenhouse gas emissions leads to a significant contribution towards global warming \cite{lashof_relative_1990}\cite{kweku_greenhouse_2018}. A straightforward solution to global warming is to reduce the level of emissions. This can be achieved by replacing fossil fuels partly with renewable energy sources, such as solar and wind energy \cite{amponsah_greenhouse_2014}\cite{kabeyi_sustainable_2022}. The downside of these sources is that it can lead to an unstable power system \cite{kiaee_utilisation_2014}\cite{ayivor_modelling_2018}\cite{tuinema_modelling_2020}\cite{samani_grid_2020}. The instability is caused by the variable energy production, due to the dependence on the weather. Fortunately, it is possible to prevent instability of the energy system by using electrolysers. An electrolyser can convert a surplus of electricity into hydrogen gas and heat. While an insufficient amount of energy is produced, the hydrogen gas produced by electrolysers can be utilised to offset the deficiency. However, utilising an electrolyser efficiently, requires a careful analysis on the placement and quantity of electrolysers in the energy network. The analysis is usually done by modelling the energy transport, also known as a load flow analysis. The difficulty in modelling the effect of an electrolyser accurately, is that you have to model multiple energy carriers. In our case, those are electricity, gas and heat, which results in an integrated energy network. \\
A way to model multi-carrier energy networks is by using a graph-based framework \cite{markensteijn_solvability_2019} \cite{markensteijn_graph-based_2020} \cite{markensteijn_convergence_2020}. The graph-based framework allows modelling a variety of coupling units, but an electrolyser has not been well researched in this framework. The closest coupling unit is a Power-to-Gas (P2G) unit and an electrical boiler, but it only captures conversion to 1 energy carrier, whilst an electrolyser outputs to different energy carriers. A combined heat and power plant (CHP) is another coupling unit that outputs to different energy carriers. Although, the input energy is different from an electrolyser, which is gas. \\
In this paper, we show a linear model for an electrolyser that can be utilised for steady-state load flow analysis. We show that boundary conditions exist that lead to a well-posed problem. Moreover, we include some numerical experiments to illustrate our approach.

%%%%%%%%%%%%%%%%%%%%%%%%%%%%%%%%%%%%%%%%%%%%%%%%%%%%%%%%%%%%%%%%%%%%%%%%%%%%%%%%%%%%%%%%%%%%%%%%%%%%

% Model

%%%%%%%%%%%%%%%%%%%%%%%%%%%%%%%%%%%%%%%%%%%%%%%%%%%%%%%%%%%%%%%%%%%%%%%%%%%%%%%%%%%%%%%%%%%%%%%%%%%%

\section{Model}
\label{section:model}
We model an electrolyser that can convert electricity into gas and residual heat, where we are interested in two cases. The first case considers an electrolyser where it converts gas and heat with a known output ratio. The second case considers a situation where the output ratio has to be determined depending on the energy transport surrounding the electrolyser. \\
Before we start with the model equations, we need a graph representation of an electrolyser. Therefore, we describe how an energy network can be transformed into a graph. A graph consists of nodes connected by links. Each node and link corresponds to an energy network element. For example, a node can represent a source and a link can represent a transmission line. Moreover, each node and link corresponds to a physical law. Explicitly, the nodes are associated with conservation laws and the links are associated with the physical model of the underlying network element. Table \ref{tab:single_carrier_model} shows an overview of the conservation laws and common models for each single-carrier. The models for each single-carrier represents the following: the electrical network is an AC system \cite{schavemaker_electrical_2017}\cite{markensteijn_graph-based_2020}, the gas network is a low-pressure system \cite{osiadacz_methods_1988}\cite{markensteijn_graph-based_2020} and the heat network is a closed-loop system \cite{markensteijn_graph-based_2020}.

\begin{table*}[pos=ht!]
    \centering
    %\small
    \caption{Models for electricity, gas and heat networks.}
    %\makebox[\textwidth][c]{
    \begin{tabular}{c|c|c|c}
        & \textbf{Network element} & \textbf{Description} &  \textbf{Model} \\ \hline
        \textbf{Electricity} & Node & Kirchhoff's law for active power & $P_i = -\sum_{j} P_{ij}$ \\
        & Node & Kirchhoff's law for reactive power & $Q_i = -\sum_{j} Q_{ij}$ \\
        & Link & Short transmission line (send, P) & $P_{ij} = g_{ij}|V_i|^2 - |V_i||V_j|\left(g_{ij} \cos{\delta_{ij}} + b_{ij} \sin{\delta_{ij}}\right)$ \\
        & Link & Transmission line (send, Q) & $Q_{ij} = -b_{ij}|V_i|^2 - |V_i||V_j|\left(g_{ij} \sin{\delta_{ij}} - b_{ij} \cos{\delta_{ij}}\right)$ \\
        & Link & Transmission line (receive, P) & $P_{ji} = g_{ij}|V_j|^2 - |V_i||V_j|\left(g_{ij} \cos{\delta_{ij}} - b_{ij} \sin{\delta_{ij}}\right)$ \\
        & Link & Transmission line (receive, Q) & $P_{ji} = -b_{ij}|V_j|^2 + |V_i||V_j|\left(g_{ij} \sin{\delta_{ij}} + b_{ij} \cos{\delta_{ij}}\right)$ \\ \hline
        \textbf{Gas} & Node & Conservation of mass & $q_i = \sum_{j} q_{ij}$ \\
        & Link & Pipe & $p_i - p_j = (C^g)^{-2} f |q_{ij}| q_{ij}$ \\ \hline
        \textbf{Heat} & Node & Conservation of mass & $m_i = \sum_{j} m_{ij}$ \\
        & Node & Conservation of energy  (supply) & $\sum_{l} m_{i, l} T^s_{i, l} = \sum_{j} m_{ij} T^s_{ij}$ \\
        & Node & Conservation of energy (return) & $\sum_{l} m_{i, l} T^r_{i, l} = \sum_{j} m_{ij} T^r_{ij}$ \\
        & Link & Pipe & $p_i - p_j = (C^h)^{-2} f |m_{ij}| m_{ij}$ \\
        & Link & Pipe heat loss (supply) & $T^s_{ji} = T^a + e^{=\frac{h \pi D L}{C_p m}}\left(T^s_{ij} - T^a\right)$ \\
        & Link & Pipe heat loss (return) & $T^r_{ij} = T^a + e^{=\frac{h \pi D L}{C_p m}}\left(T^r_{ji} - T^a\right)$ \\
        & Terminal link & Total heat power & $\Delta \varphi_{i,l} = C_p m_{i,l} (T^s_{i,l} - T^r_{i,l})$ \\ \hline
    \end{tabular}
    %}
    \label{tab:single_carrier_model}
\end{table*}

Since an electrolyser interacts with different energy-carriers, it is modelled as a node connected with links of the relevant energy carriers. We assume that these links have no energy losses. The coupling node, that represents an electrolyser, has a model that governs the energy balance given by Equation \eqref{eq:electrolyser_model}:
\begin{align}
    \label{eq:electrolyser_model}
    \eta P = \text{HHV} q + \Delta\varphi
\end{align}
where $\eta\in [0, 1]$ is the efficiency, $P$ is the active power, $q$ is the gas flow, HHV is the higher heating value of gas and $\Delta \varphi$ is the heat power. This model is based on the CHP model from \cite{savola_off-design_2005} \cite{werner_district_2013}, where we have adjusted our model based on the input and output energy. However, a CHP can generate electricity and heat in a flexible manner, whilst for an electrolyser only residual heat is available. This behaviour is reflected by Equation \eqref{eq:residual_heat_model}: 
\begin{align}
    \label{eq:residual_heat_model}
    \eta_{h} \eta P = \Delta \varphi
\end{align}
where $\eta_{h} \in [0, 1]$ is the heat efficiency. In other words, Equation \eqref{eq:residual_heat_model} tells us that a fraction of the available energy is converted to heat. Henceforth, we require two equations to model an electrolyser that converts electricity into gas and heat. \\
Our model allows the electrolyser to produce gas and heat in a flexible manner by letting $\eta_{h}$ to be unknown. Hence, our model is capable of modelling both cases we mentioned in the beginning of this section. From a mathematical point of view, we seek to know when the model is well-posed for these cases. Therefore, we show necessary conditions for when this holds in Section \ref{section:boundary_conditions}. 

%%%%%%%%%%%%%%%%%%%%%%%%%%%%%%%%%%%%%%%%%%%%%%%%%%%%%%%%%%%%%%%%%%%%%%%%%%%%%%%%%%%%%%%%%%%%%%%%%%%%

% Boundary conditions

%%%%%%%%%%%%%%%%%%%%%%%%%%%%%%%%%%%%%%%%%%%%%%%%%%%%%%%%%%%%%%%%%%%%%%%%%%%%%%%%%%%%%%%%%%%%%%%%%%%%

\section{Boundary conditions}
\label{section:boundary_conditions}
\begin{figure}[pos=ht!]
    \centering
    \begin{tikzpicture}[scale=0.75]
        %Electricity        
        \node[shape=circle, draw=orange, fill=orange!20, minimum size=5mm, inner sep=0pt, label=right:{}] (0e) at (4,0) {$0^e$};
        %Terminal links
        \path[->, draw=orange, very thick] (0e) edge node[] {} +(0,+1.5);
        %Terminal link labels
        \node[shape=circle, fill=orange, draw=orange, minimum size=1mm, inner sep=0pt, label=left:{$P_{0^e}$}, label=right:{$Q_{0^e}$}] (0e_label) at (4,0.75) {};
    
        %Gas
        \node[shape=circle, draw=gray, fill=gray!20, minimum size=5mm, inner sep=0pt, label=left:{}] (0g) at (-4,0) {$0^g$};   
        %Terminal links
        \path[->, draw=gray, very thick] (0g) edge node[] {} +(0,+1.5);
        %Terminal link labels
        \node[draw=none, label=left:{$q_{0^g}$}] (0g_label) at (-4,0.75) {};

        %Heat
        \node[shape=circle, draw=cyan, fill=cyan!20, minimum size=5mm, inner sep=0pt, label=right:{}] (0h) at (0,4) {$0^h$};
        %Terminal links
        \path[->, draw=cyan, very thick] (0h) edge node[] {} +(0,+1.5);
        %Terminal link labels
        \node[shape=circle, fill=cyan, draw=cyan, minimum size=1mm, inner sep=0pt, label=left:{$m_{0^h}, T_{0^h, l}^{s}, T_{0^h, l}^{r}$}, label=right:{$\Delta \varphi_{0^h, l}$}] (0h_label) at (0,4.75) {};

        %Coupling
        \node[shape=circle, draw=green, fill=green!20, minimum size=5mm, inner sep=0pt] (0c) at (0,0) {$0^c$};
        %Link to main nodes - electricity
        \path[-, above, draw=green, very thick] (0c) edge[] node[] {} (0e);
        %Link labels to main nodes - electricity
        \node[shape=circle, fill=green, draw=green, minimum size=1mm, inner sep=0pt, label=above:{}, label=below:{}] (0c0e) at (1,0) {};
        \node[shape=circle, fill=green, draw=green, minimum size=1mm, inner sep=0pt, label=above:{$P_{0^e 0^c}$}, label=below:{$Q_{0^e 0^c}$}] (0e0c) at (3,0) {};
        %Link to main nodes - gas
        \path[-, above, draw=green, very thick] (0c) edge[->-=.5] node[] {$q_{0^c 0^g}$} (0g);
        %Link to main nodes - heat
        \path[-, right, draw=green, very thick] (0c) edge[->-=.5] node[] {$m_{0^c 0^h}$} (0h);
        %Link labels to main nodes - heat
        \node[shape=circle, fill=green, draw=green, minimum size=1mm, inner sep=0pt, label=left:{$T_{0^c 0^h}^s, T_{0^c 0^h}^r$}, label=right:{$\Delta \varphi_{0^c 0^h}$}] (0c0h) at (0,1) {};
        \node[shape=circle, fill=green, draw=green, minimum size=1mm, inner sep=0pt, label=left:{}, label=right:{}] (0h0c) at (0,3) {};
    \end{tikzpicture}
    \caption{A graph representation of an electrolyser. Node $0^c$ and the connecting dummy links represent the electrolyser. Node $0^g$ is a node of the gas network. Node $0^h$ is a node of the heat network. Node $0^e$ is a node of the electrical network. For each node, except node $0^c$, a terminal link is connected to it. This link represents energy flowing in or out of the network.}
    \label{fig:electrolyser}
\end{figure}
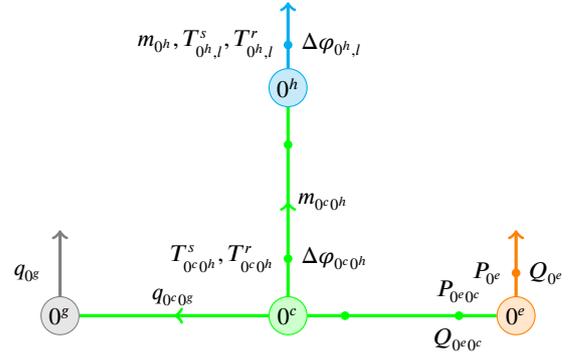

In general, integrated energy networks lead to a system of nonlinear equations. These systems are usually solved numerically. Solving a nonlinear energy system is conventionally done with the Newton-Raphson (NR) method. In this method a linearisation of the system of equations is involved, where one has to solve a linear system in order to solve the original system. With regards to linear systems, a necessary condition for well-posedness is that the system of equations is square. Usually, energy network models have more variables than equations, unless we allow additional conditions in the form of specifying variables. In this paper, we denote them as boundary conditions. A side note, nodes with different types of boundary conditions are known as node types in literature related to electrical networks. \\
To understand which conditions need to be met for a well-posed system including an electrolyser, we start with a simple network. We investigate an electrolyser with one node of each single-carrier attached to it. The graph network of this is shown in Figure \ref{fig:electrolyser}. 

\begin{figure}[pos=ht!]
    \centering
    \begin{tikzpicture}
        %Coupling
        \node[shape=circle, draw=green, fill=green!20, minimum size=5mm, inner sep=0pt] (0c) at (0,0) {$0^c$};

        %Heat
        \node[shape=circle, draw=cyan, fill=cyan!20, minimum size=5mm, inner sep=0pt, label=left:{}, label=right:{}] (0h) at (0,4) {$0^h$};
        %Link to main nodes
        \path[-, right, draw=green, very thick] (0c) edge[->-=.5] node[] {$m_{0^c 0^h}$} (0h);
        %Terminal links
        \path[->, draw=cyan, very thick] (0h) edge node[] {} +(0,+1.5);
        %Link labels
        \node[shape=circle, fill=green, draw=green, minimum size=1mm, inner sep=0pt, label=left:{$T_{0^c 0^h}^s, T_{0^c 0^h}^r$}, label=right:{$\Delta \varphi_{0^c 0^h}$}] (0c0h) at (0,1) {};
        \node[shape=circle, fill=green, draw=green, minimum size=1mm, inner sep=0pt, label=left:{}, label=right:{}] (0h0c) at (0,3) {};
        \node[shape=circle, fill=cyan, draw=cyan, minimum size=1mm, inner sep=0pt, label=left:{$m_{0^h}, T_{0^h, l}^{s}, T_{0^h, l}^{r}$}, label=right:{$\Delta \varphi_{0^h, l}$}] (0h_label) at (0,4.75) {};

        %Electricity        
        \node[shape=circle, draw=orange, fill=orange!20, minimum size=5mm, inner sep=0pt, label=right:{}] (0e) at (4,0) {$0^e$};
        %Terminal links
        \path[->, draw=orange, very thick] (0e) edge node[] {} +(0,+1.5);
        %Link labels
        \node[shape=circle, fill=green, draw=green, minimum size=1mm, inner sep=0pt, label=above:{}, label=below:{}] (0c0e) at (1,0) {};
        \node[shape=circle, fill=green, draw=green, minimum size=1mm, inner sep=0pt, label=above:{$P_{0^e 0^c}$}, label=below:{$Q_{0^e 0^c}$}] (0e0c) at (3,0) {};
        \node[shape=circle, fill=orange, draw=orange, minimum size=1mm, inner sep=0pt, label=left:{$P_{0^e}$}, label=right:{$Q_{0^e}$}] (0e_label) at (4,0.75) {};
        %Links to main nodes
        \path[-, above, draw=green, very thick] (0c) edge[] node[] {} (0e);
    \end{tikzpicture}
    \caption{A coupling between an electrical and a heat network with an electrolyser. This is equivalent to modelling an electrical boiler}
    \label{fig:electrolyser_eb}
\end{figure}
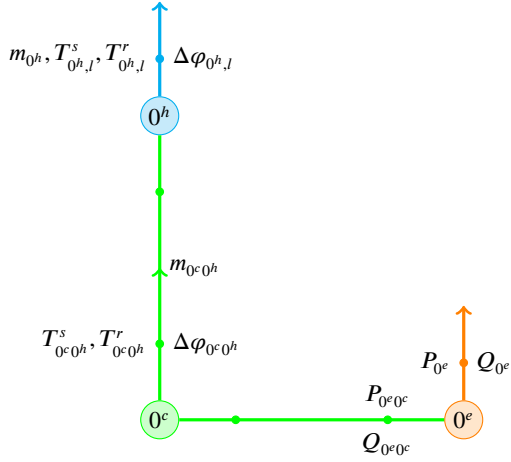

Well-posedness conditions must hold for elementary cases, such as an electrolyser only generating one of the energy outputs, gas or heat. This reduces down to a P2G unit or an electrical boiler. Leading to an analysis of two networks, which are shown in Figure \ref{fig:electrolyser_p2g} and Figure \ref{fig:electrolyser_eb}.

\begin{figure}[pos=ht!]
    \centering
    \begin{tikzpicture}[scale=0.75]
        %Coupling
        \node[shape=circle, draw=green, fill=green!20, minimum size=5mm, inner sep=0pt] (0c) at (0,0) {$0^c$};
    
        %Gas
        \node[shape=circle, draw=gray, fill=gray!20, minimum size=5mm, inner sep=0pt, label=left:{}] (0g) at (-4,0) {$0^g$};
        %Link labels
        \node[draw=none, label=left:{$q_{0^g}$}] (0g_label) at (-4,0.75) {};
        %Links to main nodes
        \path[-, above, draw=green, very thick] (0c) edge[->-=.5] node[] {$q_{0^c 0^g}$} (0g);
        %Terminal links
        \path[->, draw=gray, very thick] (0g) edge node[] {} +(0,+1.5);

        %Electricity        
        \node[shape=circle, draw=orange, fill=orange!20, minimum size=5mm, inner sep=0pt, label=right:{}] (0e) at (4,0) {$0^e$};
        %Terminal links
        \path[->, draw=orange, very thick] (0e) edge node[] {} +(0,+1.5);
        %Link labels
        \node[shape=circle, fill=green, draw=green, minimum size=1mm, inner sep=0pt, label=above:{}, label=below:{}] (0c0e) at (1,0) {};
        \node[shape=circle, fill=green, draw=green, minimum size=1mm, inner sep=0pt, label=above:{$P_{0^e 0^c}$}, label=below:{$Q_{0^e 0^c}$}] (0e0c) at (3,0) {};
        \node[shape=circle, fill=orange, draw=orange, minimum size=1mm, inner sep=0pt, label=left:{$P_{0^e}$}, label=right:{$Q_{0^e}$}] (0e_label) at (4,0.75) {};
        
        %Links to main nodes
        \path[-, above, draw=green, very thick] (0c) edge[] node[] {} (0e);
    \end{tikzpicture}
    \caption{A coupling between an electrical and a gas network with an electrolyser. This is equivalent to modelling a P2G unit.}
    \label{fig:electrolyser_p2g}
\end{figure}
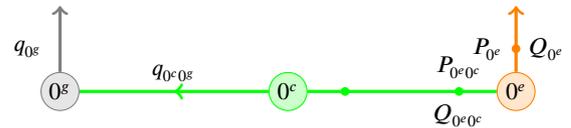

Based on our initial model, we can model a P2G unit or an electrical boiler by setting the heat efficiency to 0 or 1. Alternatively, the mass flow $m_{0^c 0^h}$ or the gas flow $q_{0^e 0^g}$ is set to $0$. \\
In the remainder of this section, we derive boundary conditions for 2 cases, a known heat efficiency and unknown heat efficiency.

%%%%%%%%%%%%%%%%%%%%%%%%%%%%%%%%%%%%%%%%%%%%%%%%%%%%%%%%%%%%%%%%%%%%%%%%%%%%%%%%%%%%%%%%%%%%%%%%%%%%

% Boundary conditions - Known heat efficiency

%%%%%%%%%%%%%%%%%%%%%%%%%%%%%%%%%%%%%%%%%%%%%%%%%%%%%%%%%%%%%%%%%%%%%%%%%%%%%%%%%%%%%%%%%%%%%%%%%%%%

\subsection{Known heat efficiency}
\label{section:known_heat_efficiency}
In this section, we assume that the heat efficiency $\eta_h$ of the electrolyser is given. Furthermore, we are interested in the production of gas and heat based on a given active power. It is also possible to specify the gas flow or heat power instead of the active power, which gives us a similar analysis. Henceforth, we focus on one case, which is the case with a specified active power. We start with deriving the boundary conditions for the P2G unit, then the electrical boiler. Finally, we combine the boundary conditions leading to well-posedness conditions for the electrolyser.

%%%%%%%%%%%%%%%%%%%%%%%%%%%%%%%%%%%%%%%%%%%%%%%%%%%%%%%%%%%%%%%%%%%%%%%%%%%%%%%%%%%%%%%%%%%%%%%%%%%%

% Boundary conditions - Known heat efficiency - P2G

%%%%%%%%%%%%%%%%%%%%%%%%%%%%%%%%%%%%%%%%%%%%%%%%%%%%%%%%%%%%%%%%%%%%%%%%%%%%%%%%%%%%%%%%%%%%%%%%%%%%

\subsubsection*{P2G}
A P2G unit has the following model \cite{clegg_integrated_2015} \cite{jiang_steady-state_2022}:
\begin{align}
    \label{eq:p2g_model}
    \eta P = \text{HHV} q
\end{align}
where $\eta \in [0, 1]$. One can obtain an equivalent model with the electrolyser model \eqref{eq:electrolyser_model}-\eqref{eq:residual_heat_model} by assuming that $\eta_{h} = 0$. The system of equations corresponding to the network shown in Figure \ref{fig:electrolyser_p2g} is represented as:
\begin{align}
    \label{eq:p2g_0}
    P_{0^e} + P_{0^e 0^c} &= 0 \\
    \label{eq:p2g_1}
    Q_{0^e} + Q_{0^e 0^c} &= 0 \\
    \label{eq:p2g_2}
    q_{0^c 0^g} - q_{0^g} &= 0 \\
    \label{eq:p2g_3}
    \eta P_{0^e 0^c} - \text{HHV} q_{0^c 0^g} &= 0
\end{align}
We require a square system for well-posedness. HHV is a known parameter, so we have 4 equations and 6 unknown variables. To obtain a square system we must specify 2 variables. \\
We assume that the input energy, $P_{0^e}$, is known. We let $Q_{0^e, 0^c} = 0$. This choice is motivated from a physical point of view. Note that we have an AC model for the electrical network. The electrolyser model \eqref{eq:p2g_model} does not depend on the reactive power, because we assume that an electrolyser requires a DC input \cite{tuinema_modelling_2020}. We also assume that there are no energy losses from converting AC to DC. Hence, the reactive power can be any arbitrary constant, we have chosen $0$ for convenience. With these assumptions, we now have a square system and the boundary conditions are presented in Table \ref{tab:electrolyser_p2g_known_heat_efficiency}. 

\begin{table}[pos=ht!]
    \centering
    \caption{Known heat efficiency: boundary conditions for a P2G unit.}
    \begin{tabular}{c|c|c} \hline
        \textbf{Node} & \textbf{Known} & \textbf{Unknown} \\ \hline
        $0^g$ & & $q_{0^g}$ \\ \hline
        $0^e$ & $P_{0^e}$ & $Q_{0^e}$ \\ \hline
        $0^c$ & $Q_{0^e 0^c}=0$ & $q_{0^c 0^g}, P_{0^c 0^e}$
    \end{tabular}
    \label{tab:electrolyser_p2g_known_heat_efficiency}
\end{table}

Note that the system is linear. With our boundary conditions, the matrix resulting from this system is non-singular. Therefore, a unique solution exists. Thus our boundary conditions lead to a well-posed problem.

%%%%%%%%%%%%%%%%%%%%%%%%%%%%%%%%%%%%%%%%%%%%%%%%%%%%%%%%%%%%%%%%%%%%%%%%%%%%%%%%%%%%%%%%%%%%%%%%%%%%

% Boundary conditions - Known heat efficiency - Electrical boiler

%%%%%%%%%%%%%%%%%%%%%%%%%%%%%%%%%%%%%%%%%%%%%%%%%%%%%%%%%%%%%%%%%%%%%%%%%%%%%%%%%%%%%%%%%%%%%%%%%%%%

\subsubsection*{Electrical Boiler}
An electrical boiler is modelled with \cite{zhao_multi-objective_2023}:
\begin{align}
    \label{eq:eb_model}
    \Delta \varphi_{0^c 0^h} &= \eta P_{0^e 0^c}
\end{align}
This model is equivalent to the electrolyser model \eqref{eq:electrolyser_model}-\eqref{eq:residual_heat_model} whenever $\eta_{h} = 1$, because this results in $q_{0^c 0^h} = 0$. We consider the model described in Equation \eqref{eq:eb_model} for the sake of brevity. The network shown in Figure \ref{fig:electrolyser_eb} results in the following system of equations:
\begin{align}
    \label{eq:eb_0}
    P_{0^e} + P_{0^e 0^c} &= 0 \\
    \label{eq:eb_1}
    Q_{0^e} + Q_{0^e 0^c} &= 0 \\
    \label{eq:eb_2}
    m_{0^c 0^h} - m_{0^h} &= 0 \\
    \label{eq:eb_3}
    m_{0^c 0^h} T^s_{0^c 0^h} - m_{0^h} T^{s}_{0^h, l} &= 0 \\
    \label{eq:eb_4}
    -m_{0^c 0^h} T^r_{0^c 0^h} + m_{0^h} T^{r}_{0^h, l} &= 0 \\
    \label{eq:eb_5}
    C_p m_{0^h} (T^{s}_{0^h, l} - T^{r}_{0^h, l}) - \Delta \varphi_{0^h, l} &= 0 \\
    \label{eq:eb_6}
    \eta P_{0^e 0^c} - \Delta \varphi_{0^c 0^h} &= 0 \\
    \label{eq:eb_7}
    C_p m_{0^c 0^h} (T_{0^c 0^h}^s - T_{0^c 0^h}^r)  - \Delta \varphi_{0^c 0^h} &= 0
\end{align}

Equation \eqref{eq:eb_7} is a new addition compared to the P2G case. This equation is required, because it describes how the heat consumption is related to the mass flow and temperature for a heat sink and source. \\
We assume that the specified heat constant $C_p$ is known. It follows that our system has $8$ equations and $12$ unknowns. To obtain a square system, we have to specify $4$ variables. 
We let the electrolyser generate heat with a specified active power $P_{0^e}$. Resulting in node $0^e$ to be a load node. From equations \eqref{eq:eb_5} and \eqref{eq:eb_7}, it follows that a reference temperature is required for a unique solution. Thus, we assume that the return temperature $T^r_{0^h, l}$ is specified.
The third variable to be specified is the supply temperature of the electrical boiler $T^s_{0^c 0^h}$, because we assume that the provided heat comes out at a set temperature. 
The last variable we specify is the reactive power,$Q_{0^e 0^c} = 0$, we do this for the same reason as for the P2G unit. The assumptions are summarised in Table \ref{tab:electrolyser_eb_known_heat_efficiency}.

\begin{table}[pos=ht!]
    \centering
    \caption{Known heat efficiency: boundary conditions for an electrical boiler.}
    \begin{tabular}{c|c|c} \hline
        \textbf{Node} & \textbf{Known} & \textbf{Unknown} \\ \hline
        $0^h$ & $T^r_{0^h, l}$ & $m_{0^h}, T^s_{0^h, l}, \Delta \varphi_{0^h, l}$ \\ \hline
        $0^e$ & $P_{0^e}$ & $Q_{0^e}$ \\ \hline
        $0^c$ & $Q_{0^c 0^e}=0, T_{0^c 0^h}^s$ & 
        \makecell{$m_{0^c 0^h}, T_{0^c 0^h}^{r}, \Delta\varphi_{0^c 0^h}, P_{0^c 0^e}$}
    \end{tabular}
    \label{tab:electrolyser_eb_known_heat_efficiency}
\end{table}

These boundary conditions lead to a unique solution. The derivation is shown below.

\begin{enumerate}
    \item The reactive power $Q_{0^e}$ is obtained from Equation \eqref{eq:eb_1}:
    \begin{align*}
        Q_{0^e} = - Q_{0^e 0^c}
    \end{align*}
    \item From Equation \eqref{eq:eb_0} we obtain the active power $P_{0^e}$:
    \begin{align*}
        P_{0^e} = -P_{0^e 0^c}
    \end{align*}
    \item Equation \eqref{eq:eb_6} yields the heat power $\Delta \varphi_{0^c 0^h}$:
    \begin{align*}
        \Delta \varphi_{0^c 0^h} = \eta P_{0^e 0^c}
    \end{align*}
    \item The supply temperature $T^s_{0^h, l}$ is obtained from Equation \eqref{eq:eb_3}:
    \begin{align*}
        T^s_{0^h, l} = \frac{m_{0^c 0^h}}{m_{0^h}} T_{0^c 0^h}^{s} \stackrel{\eqref{eq:eb_2}}{=} T_{0^c 0^h}^{s}
    \end{align*}
    \item From Equation \eqref{eq:eb_4}, we express the return temperature $T_{0^c 0^h}^r$ as:
    \begin{align*}
        T^r_{0^c 0^h} = \frac{m_{0^h}}{m_{0^c 0^h}} T^r_{0^h, l} \stackrel{\eqref{eq:eb_2}}{=} T^r_{0^h, l}
    \end{align*}
    \item Applying the results provided in steps 4 and 5 on equations \eqref{eq:eb_5} and \eqref{eq:eb_7} yields the total heat power $\Delta \varphi_{0^c 0^h}$:
    \begin{align*}
        \Delta \varphi_{0^c 0^h} = \Delta \varphi_{0^h, l}
    \end{align*}
    \item The mass flow $m_{0^c 0^h}$ is obtained from Equation \eqref{eq:eb_7}:
    \begin{align*}
        m_{0^c 0^h} = \frac{\Delta \varphi_{0^c 0^h}}{C_p (T^s_{0^c 0^h} - T^r_{0^c 0^h})}
    \end{align*}
    \item The mass flow $m_{0^h}$ is derived from Equation \eqref{eq:eb_2}:
    \begin{align*}
        m_{0^h} = m_{0^c 0^h}
    \end{align*}
\end{enumerate}

Thus all variables can be determined uniquely. Henceforth, we can conclude that the conditions shown in Table \ref{tab:electrolyser_eb_known_heat_efficiency} lead to a well-posed problem. 

%%%%%%%%%%%%%%%%%%%%%%%%%%%%%%%%%%%%%%%%%%%%%%%%%%%%%%%%%%%%%%%%%%%%%%%%%%%%%%%%%%%%%%%%%%%%%%%%%%%%

% Boundary conditions - Known heat efficiency - Electrolyser

%%%%%%%%%%%%%%%%%%%%%%%%%%%%%%%%%%%%%%%%%%%%%%%%%%%%%%%%%%%%%%%%%%%%%%%%%%%%%%%%%%%%%%%%%%%%%%%%%%%%

\subsubsection*{Electrolyser}
We have derived boundary conditions for the P2G unit and electrical boiler, which can be seen as special cases of the electrolyser with one output energy. Now we consider the electrolyser with both its output capabilities. The system of equations is given below:

\begin{align}
    \label{eq:electrolyser_0}
    P_{0^e} + P_{0^e 0^c} &= 0 \\
    \label{eq:electrolyser_1}
    Q_{0^e} + Q_{0^e 0^c} &= 0 \\
    \label{eq:electrolyser_2}
    q_{0^c 0^g} - q_{0^g} &= 0 \\
    \label{eq:electrolyser_3}
    m_{0^c 0^h} - m_{0^h} &= 0 \\
    \label{eq:electrolyser_4}
    m_{0^c 0^h} T^s_{0^c 0^h} - m_{0^h} T^{s}_{0^h, l} &= 0 \\
    \label{eq:electrolyser_5}
    -m_{0^c 0^h} T^r_{0^c 0^h} + m_{0^h} T^{r}_{0^h, l} &= 0 \\
    \label{eq:electrolyser_6}
    C_p m_{0^h} (T^{s}_{0^h, l} - T^{r}_{0^h, l}) - \Delta \varphi_{0^h, l} &= 0 \\
    \label{eq:electrolyser_7}
    C_p m_{0^c 0^h} (T_{0^c 0^h}^s - T_{0^c 0^h}^r)  - \Delta \varphi_{0^c 0^h} &= 0 \\
    \label{eq:electrolyser_8}
    \eta P_{0^e 0^c} - \text{HHV} q_{0^c 0^g} - \Delta \varphi_{0^c 0^h} &= 0 \\
    \label{eq:electrolyser_9}
    \eta_{h} \eta P_{0^e 0^c} - \Delta \varphi_{0^c 0^h} &= 0
\end{align}

Recall that $\eta_{h}$ is known. This leads to a system with 10 equations and 14 unknowns. Therefore, we have to specify 4 variables. We combine the boundary conditions derived for the P2G unit and electrical boiler to obtain boundary conditions for the electrolyser. This results in the conditions shown in Table \ref{tab:electrolyser_known_heat_efficiency}.

\begin{table}[pos=ht!]
    \centering
    \caption{Known heat efficiency: boundary conditions for an electrolyser.}
    %\makebox[\textwidth][c]{
    \begin{tabular}{c|c|c} \hline
        \textbf{Node} & \textbf{Known} & \textbf{Unknown} \\ \hline
        $0^e$ & $P_{0^e}$ & $Q_{0^e}$ \\ \hline
        $0^g$ & & $q_{0^g}$ \\ \hline
        $0^h$ & $T^r_{0^h, l}$ & \makecell{$m_{0^h}, T^s_{0^h, l}, \Delta \varphi_{0^h, l}$} \\ \hline
        $0^c$ & $Q_{0^e 0^c}=0, T^s_{0^c 0^h}$ & \makecell{$q_{0^c 0^g}, P_{0^c 0^e}, m_{0^c 0^h}, T_{0^c 0^h}^{r}, \Delta\varphi_{0^c 0^h}$}
    \end{tabular}
   % }
    \label{tab:electrolyser_known_heat_efficiency}
\end{table}

A unique solution can be obtained in a similar fashion as for the electrical boiler. We conclude that the boundary conditions lead to a well-posed problem.

%%%%%%%%%%%%%%%%%%%%%%%%%%%%%%%%%%%%%%%%%%%%%%%%%%%%%%%%%%%%%%%%%%%%%%%%%%%%%%%%%%%%%%%%%%%%%%%%%%%%

% Boundary conditions - Unknown heat efficiency

%%%%%%%%%%%%%%%%%%%%%%%%%%%%%%%%%%%%%%%%%%%%%%%%%%%%%%%%%%%%%%%%%%%%%%%%%%%%%%%%%%%%%%%%%%%%%%%%%%%%

\subsection{Unknown heat efficiency}
\label{section:unknown_heat_efficiency}
Our model allows the electrolyser to output gas and heat with any arbitrary ratio, which is equivalent to an unknown heat efficiency $\eta_{h}$. Compared to the case with a known heat efficiency, we have to specify an additional variable. Otherwise there are an infinite amount of choices for the output ratio of gas and heat.
We model a case where both output energies are known, so that we have to compute the required input energy for the electrolyser. This leads to the active power to be unknown and the gas flow and heat power to be known. We note that other cases with 2 specified energy streams can be chosen, but those lead to a similar analysis. \\
We have the same set of equations as in the case with a known heat efficiency. With an unknown heat efficiency, the system has 10 equations and 15 unknowns, so 5 variables have to be specified. Since we want the output energies to be known, we let node $0^g$ to be a load node and node $0^h$ to be a sink. Hence, we specify $q_{0^g}$ and $\Delta \varphi_{0^h, l}$. The reactive power $Q_{0^e 0^c}$, the coupling supply temperature $T^s_{0^c 0^h}$ and the return temperature $T^r_{0^h , l}$ are known. These variables follow the same reasoning as for the case with a known heat efficiency. The boundary conditions are shown in table \ref{tab:electrolyser_unknown_heat_efficiency}.

\begin{table}[pos=ht!]
    \centering
    \caption{Unknown heat efficiency: boundary conditions for an electrolyser.}
    \begin{tabular}{c|c|c} \hline
        \textbf{Node} & \textbf{Known} & \textbf{Unknown} \\ \hline
        $0^g$ & $q_{0^g}$ & \\ \hline
        $0^h$ & $\Delta \varphi_{0^h, l} > 0, T^r_{0^h, l}$ & $m_{0^h}, T^s_{0^h, l}$ \\ \hline
        $0^e$ & $P_{0^e}, Q_{0^e}$ \\ \hline
        $0^c$ & $Q_{0^e 0^c}=0, T_{0^c 0^h}^s$ & \makecell{$q_{0^c 0^g}, P_{0^c 0^e}, m_{0^c 0^h}, T_{0^c 0^h}^{r}, \Delta\varphi_{0^c 0^h}$}
    \end{tabular}
    \label{tab:electrolyser_unknown_heat_efficiency}
\end{table}

Applying these boundary conditions for our system of equations and solving it in a similar fashion as we have seen in section \ref{section:known_heat_efficiency} for the electrical boiler, we obtain a unique solution. Thus the boundary conditions results in a well-posed problem.

%%%%%%%%%%%%%%%%%%%%%%%%%%%%%%%%%%%%%%%%%%%%%%%%%%%%%%%%%%%%%%%%%%%%%%%%%%%%%%%%%%%%%%%%%%%%%%%%%%%%

% Boundary conditions - Electrolyser with physical links

%%%%%%%%%%%%%%%%%%%%%%%%%%%%%%%%%%%%%%%%%%%%%%%%%%%%%%%%%%%%%%%%%%%%%%%%%%%%%%%%%%%%%%%%%%%%%%%%%%%%

\subsection{Electrolyser coupled with single-carrier systems}
\label{section:electrolyser_coupled_with_single-carrier_systems}
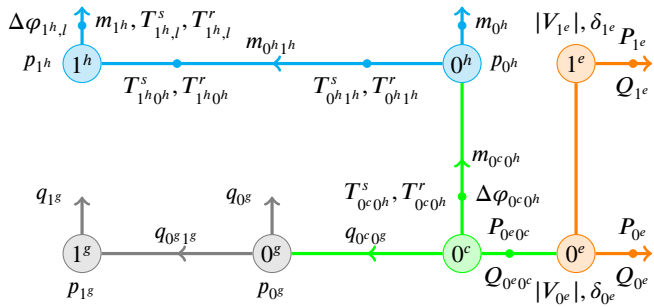
\begin{figure}[pos=ht!]
    \centering
    \begin{tikzpicture}[scale=1]
        %Electricity        
        \node[shape=circle, draw=orange, fill=orange!20, minimum size=5mm, inner sep=0pt, label=below:{$|V_{0^e}|, \delta_{0^ e}$}] (0e) at (1.5,0) {$0^e$};
        \node[shape=circle, draw=orange, fill=orange!20, minimum size=5mm, inner sep=0pt, label=above:{$|V_{1^e}|, \delta_{1^e}$}] (1e) at (1.5,2.5) {$1^e$};
        %Terminal links
        \path[->, draw=orange, very thick] (0e) edge node[] {} +(1,0);
        \path[->, draw=orange, very thick] (1e) edge node[] {} +(1,0);
        %Terminal link labels
        \node[shape=circle, fill=orange, draw=orange, minimum size=1mm, inner sep=0pt, label=above:{$P_{0^e}$}, label=below:{$Q_{0^e}$}] (0e_label) at (2.25,0) {};
        \node[shape=circle, fill=orange, draw=orange, minimum size=1mm, inner sep=0pt, label=above:{$P_{1^e}$}, label=below:{$Q_{1^e}$}] (0e_label) at (2.25,2.5) {};
        % %Link labels
        % \node[shape=circle, fill=orange, draw=orange, minimum size=1mm, inner sep=0pt, label=right:{$P_{0^e 1^e}, Q_{0^e 1^e}$}] (0h1h) at (5, 1.25) {};
        % \node[shape=circle, fill=orange, draw=orange, minimum size=1mm, inner sep=0pt, label=right:{$P_{1^e 0^e}, Q_{1^e 0^e}$}] (1h0h) at (5, 3.75) {};

        %Gas
        \node[shape=circle, draw=gray, fill=gray!20, minimum size=5mm, inner sep=0pt, label=below:{$p_{0^g}$}] (0g) at (-2.5,0) {$0^g$};
        \node[shape=circle, draw=gray, fill=gray!20, minimum size=5mm, inner sep=0pt, label=below:{$p_{1^g}$}] (1g) at (-5,0) {$1^g$};
        %Terminal links
        \path[->, draw=gray, very thick] (0g) edge node[] {} +(0,+0.75);
        \path[->, draw=gray, very thick] (1g) edge node[] {} +(0,+0.75);
        %Terminal link labels
        \node[draw=none, label=left:{$q_{0^g}$}] (0g_label) at (-2.5,0.75) {};
        \node[draw=none, label=left:{$q_{1^g}$}] (1g_label) at (-5,0.75) {};

        %Heat
        \node[shape=circle, draw=cyan, fill=cyan!20, minimum size=5mm, inner sep=0pt, label=right:{$p_{0^h}$}] (0h) at (0,2.5) {$0^h$};
        \node[shape=circle, draw=cyan, fill=cyan!20, minimum size=5mm, inner sep=0pt, label=left:{$p_{1^h}$}] (1h) at (-5,2.5) {$1^h$};
        %Terminal links
        \path[->, draw=cyan, very thick] (0h) edge node[] {} +(0,+0.75);
        \path[->, draw=cyan, very thick] (1h) edge node[] {} +(0,+0.75);
        %Terminal link labels
        \node[shape=circle, fill=cyan, draw=cyan, minimum size=1mm, inner sep=0pt, label=left:{}, label=right:{$m_{0^h}$}] (0h_label) at (0,3) {};
        \node[shape=circle, fill=cyan, draw=cyan, minimum size=1mm, inner sep=0pt, label=left:{$\Delta \varphi_{1^h, l}$}, label=right:{$m_{1^h}, T^s_{1^h, l}, T^r_{1^h, l}$}] (1h_label) at (-5,3) {};
        %Link labels
        \node[shape=circle, fill=cyan, draw=cyan, minimum size=1mm, inner sep=0pt, label=below:{$T_{0^h 1^h}^s, T_{0^h 1^h}^r$}] (0h1h) at (-1.25, 2.5) {};
        \node[shape=circle, fill=cyan, draw=cyan, minimum size=1mm, inner sep=0pt, label=below:{$T_{1^h 0^h}^s, T_{1^h 0^h}^r$}] (1h0h) at (-3.75, 2.5) {};
        
        %Coupling
        \node[shape=circle, draw=green, fill=green!20, minimum size=5mm, inner sep=0pt] (0c) at (0,0) {$0^c$};
        %Link to main nodes - electricity
        \path[-, above, draw=green, very thick] (0c) edge[] node[] {} (0e);
        \path[-, above, draw=orange, very thick] (0e) edge[] node[] {} (1e);
        %Link labels to main nodes - electricity
%        \node[shape=circle, fill=green, draw=green, minimum size=1mm, inner sep=0pt, label=above:{}, label=below:{}] (0c0e) at (0.625,0) {};
        \node[shape=circle, fill=green, draw=green, minimum size=1mm, inner sep=0pt, label=above:{$P_{0^e 0^c}$}, label=below:{$Q_{0^e 0^c}$}] (0e0c) at (0.625,0) {};
        %Link to main nodes - gas
        \path[-, above, draw=green, very thick] (0c) edge[->-=.5] node[] {$q_{0^c 0^g}$} (0g);
        \path[-, above, draw=gray, very thick] (0g) edge[->-=.5] node[] {$q_{0^g 1^g}$} (1g);
        %Link to main nodes - heat
        \path[-, right, draw=green, very thick] (0c) edge[->-=.5] node[] {$m_{0^c 0^h}$} (0h);
        \path[-, above, draw=cyan, very thick] (0h) edge[->-=.5] node[] {$m_{0^h 1^h}$} (1h);
        %Link labels to main nodes - heat
        \node[shape=circle, fill=green, draw=green, minimum size=1mm, inner sep=0pt, label=left:{$T_{0^c 0^h}^s, T_{0^c 0^h}^r$}, label=right:{$\Delta \varphi_{0^c 0^h}$}] (0c0h) at (0,0.75) {};
%        \node[shape=circle, fill=green, draw=green, minimum size=1mm, inner sep=0pt, label=left:{}, label=right:{}] (0h0c) at (0,1.875) {};
    \end{tikzpicture}
    \caption{The network shown in Figure \ref{fig:electrolyser} is extended with physical links for each energy carrier. Node $0^e$, $0^g$ and $0^h$ are acting as junctions. Whilst node $1^e$, $1^g$ and $1^h$ are sources or sinks.}
    \label{fig:electrolyser_physical_links}
\end{figure}

In the previous sections, we have investigated a network with one electrolyser without any physical links connected to it, so energy losses from the connected single-carrier networks are not modelled. To show the effect of an electrolyser in a more realistic energy network setting, we simply extend the network shown in Figure \ref{fig:electrolyser} with one physical link for each single-carrier network. The extended link of the gas and heat network represents a pipe. The link for the electrical network represents a transmission line. The network is shown in Figure \ref{fig:electrolyser_physical_links}.

Nodes $0^e$, $0^g$ and $0^h$ are junctions. The junctions are modelled in a conventional way by assuming that no energy can enter or escape the network. This forces our load nodes from our previous network ($0^e$, $0^g$ and $0^h$) to move to node $1^e$, $1^g$ and $1^h$ in the current network. The resulting system of equations is shown in Appendix \ref{appendix:system_of_equations}, which has 18 equations and 30 unknown variables whenever the heat efficiency is known.

For these systems, the same well-posedness conditions derived in Section \ref{section:known_heat_efficiency} with a known heat efficiency and Section \ref{section:unknown_heat_efficiency} with an unknown heat efficiency are applied for loads $1^e$, $1^g$ and $1^h$. The gas and heat network require reference pressures. The reference pressure can be placed on a junction or a load node in their respective network. We have chosen to place them on the loads. Also the coupling node $0^c$ has the same boundary conditions that we have derived before. The resulting boundary conditions that lead to well-posedness for a known and unknown heat efficiency are shown in Table \ref{tab:electrolyser_physical_links_known_heat_efficiency} and \ref{tab:electrolyser_physical_links_unknown_heat_efficiency}. 

\begin{table}[pos=ht!]
    \centering
    \caption{Known heat efficiency: boundary conditions for an electrolyser with physical links.}
    %\makebox[\textwidth][c]{
    \begin{tabular}{c|c|c} \hline
        \textbf{Node} & \textbf{Known} & \textbf{Unknown} \\ \hline
        $0^e$ & $P_{0^e}=0, Q_{0^e}=0$ & $V_{0^e}, \delta_{0^e}$ \\ \hline
        $1^e$ & $P_{1^e}, V_{1^e}, \delta_{1^e}$ & $Q_{1^e}$ \\ \hline
        $0^g$ & $q_{0^g} = 0$ & $p_{0^g}$ \\ \hline
        $1^g$ & $p_{1^g}$ & $q_{1^g}$ \\ \hline
        $0^h$ & $m_{0^h}=0$ & $p_{0^h}$ \\ \hline
        $1^h$ & $p_{1^h}, T^{r}_{1^h, l}$ & \makecell{$p_{1^h}, T^s_{1^h, l}, \Delta\varphi_{1^h, l}$} \\ \hline
        $0^c$ & $Q_{0^e 0^c}=0, T_{0^c 0^h}^s$ & $q_{0^c 0^g}, P_{0^e 0^c}, m_{0^c 0^h}, T_{0^c 0^h}^r, \Delta\varphi_{0^c 0^h}$
    \end{tabular}
   % }
    \label{tab:electrolyser_physical_links_known_heat_efficiency}
\end{table}

\begin{table}[pos=ht!]
    \centering
    \caption{Unknown heat efficiency: boundary conditions for an electrolyser with physical links.}
    \begin{tabular}{c|c|c} \hline
        \textbf{Node} & \textbf{Known} & \textbf{Unknown} \\ \hline
        $0^e$ & $P_{0^e}=0, Q_{0^e}=0$ & $V_{0^e}, \delta_{0^e}$ \\ \hline
        $1^e$ & $V_{1^e}, \delta_{1^e}$ & $P_{1^e}, Q_{1^e}$ \\ \hline
        $0^g$ & $q_{0^g} = 0$ & $p_{0^g}$  \\ \hline
        $1^g$ & $q_{1^g} > 0, p_{1^g}$ & \\ \hline
        $0^h$ & $m_{0^h}=0$ & $p_{0^h}$ \\ \hline
        $1^h$ & $p_{1^h}, \Delta\varphi_{1^h, l} > 0, T^{r}_{1^h, l}$ & \makecell{$m_{1^h}, T^s_{1^h, l}$} \\ \hline
        $0^c$ & $Q_{0^e 0^c}=0, T_{0^c 0^h}^s$ & $q_{0^c 0^g}, P_{0^e 0^c}, m_{0^c 0^h}, T_{0^c 0^h}^r, \Delta\varphi_{0^c 0^h}$
    \end{tabular}
    \label{tab:electrolyser_physical_links_unknown_heat_efficiency}
\end{table}

An analytical solution is not easy to derive, due to the nonlinear equations in our system. Therefore, we validate our results numerically.

%%%%%%%%%%%%%%%%%%%%%%%%%%%%%%%%%%%%%%%%%%%%%%%%%%%%%%%%%%%%%%%%%%%%%%%%%%%%%%%%%%%%%%%%%%%%%%%%%%%%

% Numerical results

%%%%%%%%%%%%%%%%%%%%%%%%%%%%%%%%%%%%%%%%%%%%%%%%%%%%%%%%%%%%%%%%%%%%%%%%%%%%%%%%%%%%%%%%%%%%%%%%%%%%

\section{Numerical Results}
\label{section:numerical_results}
In this section, we numerically solve the system representing an electrolyser with physical links shown in Figure \ref{fig:electrolyser_physical_links}. We do this with a known and unknown heat efficiency. Recall that there are some assumptions that hold for this energy network. Firstly, the electrical network is modelled as an AC system. Secondly, we assume a low-pressure system for the gas network. Lastly, the heat network is modelled as a closed-loop system with a supply and return line. Moreover, specification of the transmission lines, gas and pipes are summarised in Appendix \ref{appendix:system_specification}. Additionally, the efficiencies of the electrolyser are found in the same section. \\
We solve the system with the Newton-Raphson (NR) method, where the stopping criterion is defined as:
\begin{align*}
    \|F\|_2 \leq 10^{-6}
\end{align*}
For the inner solve we have used an LU factorisation from the SuperLU package \cite{demmel_superlu_1999}. \\
The system of equations has 18 equations and variables. By substitution, we are able to reduce the system to 15 equations and variables. The initial guess for the NR method is given in Table \ref{tab:initial_guess_known_heat_efficiency_pl}. We motivate our choices as follows. The mass flows are chosen such that they are nonzero, otherwise the Jacobian is singular. The supply and return temperatures are determined such that the average of these temperatures equals the temperature of the boundary condition for the supply temperature $T^s_{0^c 0^h}$. The pressure in the gas network $p_{0^g}$ is chosen such that the pressure drop is nonzero, because a pressure drop of zero leads to a singular Jacobian. Likewise for the pressure $p_{0^h}$ in the heat network. The voltage magnitude $|V_{0^e}|$ is based on a flat start with same value as the boundary condition. The other variables are set to 0 out of convenience.

\begin{table}[pos=ht!]
    \centering
    \caption{Known heat efficiency: the initial guess.}
    \begin{tabular}{ccccc}
         \textbf{Electricity} & & \textbf{Gas} \\ \hline\hline
         $|V_{0^e}|$ & $\delta_{0^e}$ & $p_{0^g}$ & $q_{0^g 1^g}$ \\ \hline
         $\frac{690}{\sqrt{3}}$V & $0$rad & $1.05$bar & $0\frac{\text{kg}}{\text{s}}$ \\ \\
         \textbf{Heat} \\ \hline\hline
         $T^{s}_{0^h 1^h}$ & $T^r_{1^h 0^h}$ & $p_{0^h}$ & $T^s_{1^h, l}$ & $m_{0^h 1^h}$ \\ \hline
         $353.15$K & $313.15$K & $6.3$bar & $353.15$K & $1\frac{\text{kg}}{\text{s}}$ \\ \\
         \textbf{Coupling} & \\ \hline\hline
         $P_{0^e 0^c}$ & $q_{0^c 0^h}$ & $m_{0^c 0^h}$ & $T^r_{0^c 0^h}$ & $\Delta\varphi_{0^c 0^h}$ \\ \hline
         $0$MW & $0\frac{\text{kg}}{\text{s}}$ & $1\frac{\text{kg}}{\text{s}}$ & $313.15$K & $0$MW
    \end{tabular}
    \label{tab:initial_guess_known_heat_efficiency_pl}
\end{table}

With this initial guess, the NR method converges in 5 iterations. The numerical solution is shown in Table \ref{tab:numerical_solution_known_heat_efficiency_pl}. 

\begin{table}[pos=ht!]
    \centering
    \caption{Known heat efficiency: solution of an electrolyser with physical links. The values of the boundary conditions are in bold.}
    \begin{tabular}{ccccc}
         \textbf{Electricity} & \\ \hline\hline
         $|V_{0^e}|$ & $\delta_{0^e}$ & $\mathbf{P_{0^e}}$ & $\mathbf{Q_{0^e}}$ & $\mathbf{|V_{1^e}|}$ \\
         $375$V & $-0.259$rad & $\mathbf{0}$\textbf{MW} & $\mathbf{0}$\textbf{MW} & $\mathbf{398}$\textbf{V} \\ \hline
         $\mathbf{\delta_{1^e}}$ & $\mathbf{P_{1^e}}$ & $Q_{1^e}$ & $P_{0^e 1^e}$ & $Q_{0^e 1^e}$ \\
         $\mathbf{0}$\textbf{rad} & $\mathbf{-2.5}$\textbf{MW} & $-0.662$MW & $-2.434$MW & $0$MW \\ \hline
         $P_{1^e 0^e}$ & $Q_{1^e 0^e}$ \\
         $2.5$MW & $0.662$MW \\ \\
         \textbf{Gas} & \\ \hline\hline
         $p_{0^g}$ & $\mathbf{q_{0^g}}$ & $\mathbf{p_{1^g}}$ & $q_{1^g}$ & $q_{0^g 1^g}$ \\
         $1.003$bar & $\mathbf{0\frac{\text{kg}}{\text{s}}}$ & $\mathbf{1}$\textbf{bar} & $0.013\frac{\text{kg}}{\text{s}}$ & $0.013\frac{\text{kg}}{\text{s}}$ \\ \\
         \textbf{Heat} & \\ \hline\hline
         $p_{0^h}$ & $\mathbf{m_{0^h}}$ & $T^s_{0^h 1^h} (T^{s}_{0^h})$ & $\mathbf{p_{1^h}}$ & $m_{1^h}$ \\ 
         $6.048$bar & $\mathbf{0}\frac{\text{kg}}{\text{s}}$ & $338.15$K & $\mathbf{6}$\textbf{bar} & $5.74\frac{\text{kg}}{\text{s}}$ \\ \hline
         $T^s_{1^h, l} (T^{s}_{1^h})$ & $\mathbf{T^{r}_{1^h, l}}$ & $T^r_{1^h 0^h} (T^{r}_{1^h})$ & $\Delta\varphi_{1^h}$ & $m_{0^h 1^h}$ \\
         $337.88$K & $\mathbf{323.15}$\textbf{K} & $323.15$K & $0.354$MW & $5.74\frac{\text{kg}}{\text{s}}$ \\ \\
         \textbf{Coupling} & \\ \hline\hline
         $P_{0^e 0^c}$ & $\mathbf{Q_{0^e 0^c}}$ & $q_{0^c 0^h}$ & $m_{0^c 0^h}$ & $\Delta\varphi_{0^c 0^h}$ \\
         $2.434$MW & $\mathbf{0}$\textbf{MW} & $0.013\frac{\text{kg}}{\text{s}}$ & $5.74\frac{\text{kg}}{\text{s}}$ & $0.365$MW \\ \hline
         $\mathbf{T^s_{0^c 0^h}}$ & $T^r_{0^c 0^h} (T^{r}_{0^h})$ \\
         $\mathbf{338.15}$\textbf{K} & $322.942$K
    \end{tabular}
    \label{tab:numerical_solution_known_heat_efficiency_pl}
\end{table}

We observe reasonable energy losses caused by the physical links. In the electrical network, the transmission line shows a loss in active power. In absolute value it drops by $0.066$MW. For the gas network, we observe a small pressure drop of $\Delta p = 0.003$bar. Similarly, for the heat network, a pressure drop of $\Delta p = 0.048$bar is noted. For the supply temperature from node $0^h$ to $1^h$, the temperature drops by $0.27$K. The return temperature in the direction from node $1^h$ to $0^h$ drops by $0.208$K. Thus our model behaves as expected. \\
For the case with an unknown heat efficiency, we have chosen boundary condition values, such that we have the same numerical solution as for the known heat efficiency case. Moreover, we have observed similar convergence behaviour with the NR method. \\
The results for the electrolyser with physical links suggest that the problem is well-posed.

%%%%%%%%%%%%%%%%%%%%%%%%%%%%%%%%%%%%%%%%%%%%%%%%%%%%%%%%%%%%%%%%%%%%%%%%%%%%%%%%%%%%%%%%%%%%%%%%%%%%

% Conclusion

%%%%%%%%%%%%%%%%%%%%%%%%%%%%%%%%%%%%%%%%%%%%%%%%%%%%%%%%%%%%%%%%%%%%%%%%%%%%%%%%%%%%%%%%%%%%%%%%%%%%

\section{Conclusion}

We have introduced a linear model for the electrolyser. This model is used within a steady-state load flow analysis of multi-carrier energy networks. Furthermore, we have focused on energy networks that consist of 3 different energy carriers. These energy carriers are electricity, gas and heat. In this context, we are interested when the inclusion of the electrolyser model lead to a well-posed problem. We have validated our results analytically as well as numerically. \\
The results of the electrolyser connected with physical links suggest that the links can be replaced with networks, since the boundary conditions connected with the electrolyser are the determining factor for well-posedness. In other words, if these boundary conditions are chosen such that they coincide with a known or unknown heat efficiency case, then well-posedness is expected for a broader set of network topologies. Therefore, our electrolyser model can be used in a larger setting than just one physical link per energy carrier.

%%%%%%%%%%%%%%%%%%%%%%%%%%%%%%%%%%%%%%%%%%%%%%%%%%%%%%%%%%%%%%%%%%%%%%%%%%%%%%%%%%%%%%%%%%%%%%%%%%%%

% Appendix

%%%%%%%%%%%%%%%%%%%%%%%%%%%%%%%%%%%%%%%%%%%%%%%%%%%%%%%%%%%%%%%%%%%%%%%%%%%%%%%%%%%%%%%%%%%%%%%%%%%%

 \appendix
%%%%%%%%%%%%%%%%%%%%%%%%%%%%%%%%%%%%%%%%%%%%%%%%%%%%%%%%%%%%%%%%%%%%%%%%%%%%%%%%%%%%%%%%%%%%%%%%%%%%

% Nomenclature

%%%%%%%%%%%%%%%%%%%%%%%%%%%%%%%%%%%%%%%%%%%%%%%%%%%%%%%%%%%%%%%%%%%%%%%%%%%%%%%%%%%%%%%%%%%%%%%%%%%%

\section{Nomenclature}
\label{appendix:nomenclature}

\begin{table}[pos=H]
    \centering
    \begin{tabular}{c|c} \hline
         \textbf{Electricity} \\ \hline\hline
         $\delta$ & Voltage angle [rad] \\
         $P$ &  Active power [W] \\
         $Q$ & Reactive power [var] \\
         $V$ & Voltage phasor [V] \\
         $|V|$ & Voltage amplitude [V] \\ 
         \hline\hline
         \textbf{Gas} \\ \hline\hline
         $f$ & Friction factor \\
         $q$ & Gas flow rate [$\text{kg} \cdot \text{s}^{-1}$] \\
         $C_g$ & Pipe constant [$\text{kg}^2 \cdot \text{m}^2$] \\
         \hline\hline
         \textbf{Heat} \\ \hline\hline
         $\varphi$ & Heat power [W] \\
         $\Delta \varphi$ & Total heat power[W] \\
         $T$ & Temperature [K] \\
         $m$ & Mass flow rate [$\text{kg} \cdot \text{s}^{-1}$] \\
         $C_h$ & Pipe constant [$\text{kg} \cdot \text{m}$] \\
         $C_p$ & Specific heat [$\text{m}^2 \cdot \text{K}^{-1} \cdot \text{s}^{-2}$] \\
         \hline\hline
         \textbf{General} \\ \hline\hline
         $p$ & Pressure [Pa] \\
    \end{tabular}
    \label{tab:list_of_symbols}
\end{table}

%%%%%%%%%%%%%%%%%%%%%%%%%%%%%%%%%%%%%%%%%%%%%%%%%%%%%%%%%%%%%%%%%%%%%%%%%%%%%%%%%%%%%%%%%%%%%%%%%%%%

% System specification

%%%%%%%%%%%%%%%%%%%%%%%%%%%%%%%%%%%%%%%%%%%%%%%%%%%%%%%%%%%%%%%%%%%%%%%%%%%%%%%%%%%%%%%%%%%%%%%%%%%%

\section{System specification}
\label{appendix:system_specification}
\begin{table}[pos=H]
    \centering
    \caption{Physical properties of the electrical system.}
    \begin{tabular}{c|c|c}
        & \textbf{Variable} & \textbf{Value} \\ \hline\hline
        & Power system & AC \\ \hline\hline
        \textbf{Line} & B (Susceptance) & $-0.3$ S \\ \hline
        & G (Conductance) & $0.03$ S \\
    \end{tabular}
    \label{tab:electricity}
\end{table}

\begin{table}[pos=H]
    \centering
    \caption{Physical properties of the gas system.}
    \begin{tabular}{c|c|c}
        & \textbf{Variable} & \textbf{Value} \\ \hline\hline
        & Pressure system & Low pressure \\ \hline
        & Gas type & Hydrogen gas \\ \hline
        & HHV & $1.418 \cdot 10^8 \frac{\text{J}}{\text{kg}}$ \\ \hline
        & S (Specific gravity) & $0.589$ \\ \hline
        & Z (Compressibility factor) & $1$ \\ \hline
        & $p_n$ (Standard pressure) & $1$ bar \\ \hline
        & $T_n$ (Standard temperature) & $288$ K \\ \hline
        & R (Ideal gas constant) & $8.314413 \frac{\text{J}}{\text{molK}}$ \\ \hline
        & M (Molar mass of air) & $28.97\cdot 10^{-3} \frac{\text{kg}}{\text{mol}}$ \\ \hline\hline
        \textbf{Pipe} & L & $500$ m \\ \hline
        & D & $0.15$ m \\ \hline
        & f (Friction factor) & $6.5 \cdot 10^{-3}$
    \end{tabular}
    \label{tab:gas}
\end{table}

\begin{table}[pos=H]
    \centering
    \caption{Physical properties of the heat system.}
    \begin{tabular}{c|c|c}
        & \textbf{Variable} & \textbf{Value} \\ \hline\hline
        & $\rho$ (density of water) & $960 \frac{\text{kg}}{\text{m$^3$}}$ \\ \hline
        & $C_p$ (Specific heat of water) & $4.182 \cdot 10^3 \frac{\text{J}}{\text{kgK}}$ \\ \hline
        & g (gravitational constant) & $9.81 \frac{\text{kg}}{\text{s$^2$}}$ \\ \hline
        & $T_a$ (ambient temperature) & $273.15$ K \\ \hline\hline
        \textbf{Pipe} & L & $500$ m \\ \hline
        & D & $0.15$ m \\ \hline
        & $\lambda$ (Heat transfer coefficient) & $0.2 \frac{\text{W}}{\text{m$^2$ K}}$ \\ \hline
        & f (Friction factor) & $6.5 \cdot 10^{-3}$
    \end{tabular}
    \label{tab:heat}
\end{table}

\begin{table}[pos=H]
    \centering
    \caption{Electrolyser efficiency.}
    \begin{tabular}{c|c}
        \textbf{Variable} & \textbf{Value} \\ \hline\hline
        $\eta$ & $\frac{9}{10}$ \\ \hline
        $\eta_{h}$ & $\frac{1}{6}$
    \end{tabular}
    \label{tab:coupling}
\end{table}

%%%%%%%%%%%%%%%%%%%%%%%%%%%%%%%%%%%%%%%%%%%%%%%%%%%%%%%%%%%%%%%%%%%%%%%%%%%%%%%%%%%%%%%%%%%%%%%%%%%%

% System of equations

%%%%%%%%%%%%%%%%%%%%%%%%%%%%%%%%%%%%%%%%%%%%%%%%%%%%%%%%%%%%%%%%%%%%%%%%%%%%%%%%%%%%%%%%%%%%%%%%%%%%

\section{System of equations of an electrolyser with physical links}
\label{appendix:system_of_equations}
The system of equations is shown below:
\begin{align}
    P_{0^e} + P_{0^e 1^e} + P_{0^e 0^c} &= 0 \\
    Q_{0^e} + Q_{0^e 1^e} + Q_{0^e 0^c} &= 0 \\
    \label{eq:kirchhoff_active_power}
    P_{1^e} + P_{1^e 0^e} &= 0 \\
    \label{eq:derived0}
    Q_{1^e} + Q_{1^e 0^e} &= 0 \\
    q_{0^c 0^g} - q_{0^g 1^g} -  q_{0^g} &= 0 \\
    \label{eq:derived1}
    q_{0^g 1^g} - q_{1^g} &= 0 \\
    p_{0^g} - p_{1^g} - \left(C^g\right)^{-2} f^g |q_{0^g 1^g}| q_{0^g 1^g} &= 0 \\
    m_{0^h 1^h} - m_{1^h} &= 0 \\
    m_{0^c 0^h} - m_{0^h 1^h} - m_{0^h} &= 0 \\
    p_{0^h} - p_{1^h} - \left(C^h\right)^{-2} f^h |m_{0^h 1^h}| m_{0^h 1^h} &= 0 \\
    m_{0^c 0^h} T^s_{0^c 0^h} - m_{0^h 1^h} T^s_{0^h 1^h} &= 0 \\
    -m_{0^c 0^h} T^r_{0^c 0^h} + m_{0^h 1^h} T^r_{0^h 1^h} &= 0 \\
    m_{0^h 1^h} T^s_{1^h 0^h} - m_{1^h} T^{s}_{1^h, l} &= 0 \\
    -m_{0^h 1^h} T^r_{1^h 0^h} + m_{1^h} T^{r}_{1^h, l} &= 0 \\
    \label{eq:derived2}
    C_p m_{1^h} (T^{s}_{1^h, l} - T^{r}_{1^h, l}) - \Delta \varphi_{1^h, l} &= 0 \\
    C_p m_{0^c 0^h} (T_{0^c 0^h}^s - T_{0^c 0^h}^r)  - \Delta \varphi_{0^c 0^h} &= 0 \\
    \eta P_{0^e 0^c} - \text{HHV} q_{0^c 0^g} - \Delta \varphi_{0^c 0^h} &= 0 \\
    \eta_h \eta P_{0^e 0^c} - \Delta \varphi_{0^c 0^h} &= 0
\end{align}
Let $\delta_{ij} = \delta_i - \delta_j$. We use $0$ and $1$ as a shorthand notation for node $0^e$ and $1^e$. We define the active and reactive power on the transmission line as:
\begin{small}
    \begin{align*}
        &P_{1 0} = g_{1 0}|V_{1}|^2 - |V_{0}||V_{1}|\left(g_{1 0} \cos{\delta_{1 0}} + b_{1 0} \sin{\delta_{1 0}}\right) \\
        &P_{0 1} = g_{0 1}|V_{0}|^2 - |V_{0}||V_{1}|\left(g_{0 1} \cos{\delta_{0 1}} - b_{0 1} \sin{\delta_{0 1}}\right) \\
        &Q_{1 0} = -b_{1 0}|V_{1}|^2 - |V_{0}||V_{1}|\left(g_{1 0} \sin{\delta_{1 0}} - b_{1 0} \cos{\delta_{1 0}}\right) \\
        &Q_{0 1} = -b_{0 1}|V_{0}|^2 + |V_{0}||V_{1}|\left(g_{0 1} \sin{\delta_{0 1}} + b_{0 1} \cos{\delta_{0 1}}\right)
    \end{align*}
\end{small}
% \begin{align*}
%     &P_{1^e 0^e} = g_{1^e 0^e}|V_{1^e}|^2 - |V_{0^e}||V_{1^e}|\left(g_{1^e 0^e} \cos{\delta_{1^e 0^e}} + b_{1^e 0^e} \sin{\delta_{1^e 0^e}}\right) \\
%     &P_{0^e 1^e} = g_{0^e 1^e}|V_{0^e}|^2 - |V_{0^e}||V_{1^e}|\left(g_{0^e 1^e} \cos{\delta_{0^e 1^e}} - b_{0^e 1^e} \sin{\delta_{0^e 1^e}}\right) \\
%     &Q_{1^e 0^e} = -b_{1^e 0^e}|V_{1^e}|^2 - |V_{0^e}||V_{1^e}|\left(g_{1^e 0^e} \sin{\delta_{1^e 0^e}} - b_{1^e 0^e} \cos{\delta_{1^e 0^e}}\right) \\
%     &Q_{0^e 1^e} = -b_{0^e 1^e}|V_{0^e}|^2 + |V_{0^e}||V_{1^e}|\left(g_{0^e 1^e} \sin{\delta_{0^e 1^e}} + b_{0^e 1^e} \cos{\delta_{0^e 1^e}}\right)
% \end{align*}
These are substituted in equations \eqref{eq:kirchhoff_active_power}-\eqref{eq:derived0}. Hence, the active and reactive power corresponding to the transmission line are no longer present in the system of equations. Instead the voltage magnitude and voltage angle are introduced into the system of equations. \\
For the heat network, assuming that $m_{0^h 1^h} > 0$, the temperature at the end of a supply line and return line are substituted in the relevant equations by:
\begin{align*}
    T^s_{1^h 0^h} &= (T^s_{0^h 1^h} - T^a)e^{\frac{-\lambda}{C_p m_{0^h 1^h}}L} + T^a \\
    T^r_{0^h 1^h} &= (T^r_{1^h 0^h} - T^a)e^{\frac{-\lambda}{C_p m_{0^h 1^h}}L} + T^a
\end{align*}
This leads to a system of 18 equations and 30 unknown variables whenever $\eta_h$ is known.

%%%%%%%%%%%%%%%%%%%%%%%%%%%%%%%%%%%%%%%%%%%%%%%%%%%%%%%%%%%%%%%%%%%%%%%%%%%%%%%%%%%%%%%%%%%%%%%%%%%%

% Credits

%%%%%%%%%%%%%%%%%%%%%%%%%%%%%%%%%%%%%%%%%%%%%%%%%%%%%%%%%%%%%%%%%%%%%%%%%%%%%%%%%%%%%%%%%%%%%%%%%%%%

\printcredits

%%%%%%%%%%%%%%%%%%%%%%%%%%%%%%%%%%%%%%%%%%%%%%%%%%%%%%%%%%%%%%%%%%%%%%%%%%%%%%%%%%%%%%%%%%%%%%%%%%%%

% Declaration of competing interest

%%%%%%%%%%%%%%%%%%%%%%%%%%%%%%%%%%%%%%%%%%%%%%%%%%%%%%%%%%%%%%%%%%%%%%%%%%%%%%%%%%%%%%%%%%%%%%%%%%%%

\section*{Declaration of competing interest}
The authors declare that they have no known competing financial interests or personal relationships that could have appeared to influence the work reported in this paper.

%%%%%%%%%%%%%%%%%%%%%%%%%%%%%%%%%%%%%%%%%%%%%%%%%%%%%%%%%%%%%%%%%%%%%%%%%%%%%%%%%%%%%%%%%%%%%%%%%%%%

% Acknowledgments

%%%%%%%%%%%%%%%%%%%%%%%%%%%%%%%%%%%%%%%%%%%%%%%%%%%%%%%%%%%%%%%%%%%%%%%%%%%%%%%%%%%%%%%%%%%%%%%%%%%%

\section*{Acknowledgments}
This work is part of the Energy Intranets (NEAT: ESI-BiDa 647.003.002) project, which is funded by the Dutch Research Council NWO, Netherlands in the framework of the Energy Systems Integration \& Big Data programme.

%%%%%%%%%%%%%%%%%%%%%%%%%%%%%%%%%%%%%%%%%%%%%%%%%%%%%%%%%%%%%%%%%%%%%%%%%%%%%%%%%%%%%%%%%%%%%%%%%%%%

% Bibliography

%%%%%%%%%%%%%%%%%%%%%%%%%%%%%%%%%%%%%%%%%%%%%%%%%%%%%%%%%%%%%%%%%%%%%%%%%%%%%%%%%%%%%%%%%%%%%%%%%%%%

\bibliographystyle{cas-model2-names}

\bibliography{main.bib}

%%%%%%%%%%%%%%%%%%%%%%%%%%%%%%%%%%%%%%%%%%%%%%%%%%%%%%%%%%%%%%%%%%%%%%%%%%%%%%%%%%%%%%%%%%%%%%%%%%%%


%%%%%%%%%%%%%%%%%%%%%%%%%%%%%%%%%%%%%%%%%%%%%%%%%%%%%%%%%%%%%%%%%%%%%%%%%%%%%%%%%%%%%%%%%%%%%%%%%%%%



\end{document}